\documentclass[12pt]{amsart}

\title{Trace identities for the Topological Vertex}
\author{Jim Bryan, Martijn Kool, Benjamin Young}
\date{\today}
\address{
\texttt{jbryan@math.ubc.ca}\\
Department of Mathematics\\
University of British Columbia \\
Room 121, 1984 Mathematics Road  \\
Vancouver, B.C., Canada V6T 1Z2  
}
\address{
\texttt{m.kool1@uu.nl}\\
Mathematical Institute\\
Utrecht University\\
Budapestlaan 6\\
3584 CD, Utrecht,
The Netherlands
}
\address{
\texttt{bjy@uoregon.edu}\\
Department of Mathematics\\
Fenton Hall\\
University of Oregon\\
Eugene, OR 97403-1222 USA
}

\usepackage{verbatim}
\usepackage{amsmath}
\usepackage{amsmath,amsthm,amsfonts,amssymb}
\usepackage{times}

\newcommand{\cnums} {{\mathbb C}}          % complex numbers
\newcommand{\znums} {{\mathbb Z}}		% integers
\newcommand{\Vsf}{\mathsf{V}}
\newcommand{\bx}{\square}
\renewcommand{\emptyset}{\varnothing}
\newcommand{\half}{\frac{1}{2}}

\newtheorem{theorem}{Theorem}%[section]

\newtheorem{lemma}[theorem]{Lemma}

\theoremstyle{definition}

\newtheorem{def-thm}[thm]{Definition-Theorem}
\newtheorem{remark}[theorem]{Remark}
\newtheorem{defn}[theorem]{Definition}

\newcommand{\tr}{\operatorname{tr}}

\newcommand{\FockSpace}{\Lambda^{\frac{\infty}{2}}V}
\newcommand{\FockSpaceZero}{\Lambda^{\frac{\infty}{2}}_{0}V}
\newcommand{\ZplusHalf}{\znums+{ \half}}
\newcommand{\E}{\mathcal{E}}
\newcommand{\ptotheminusrho}{p^{-\rho}}
\renewcommand{\O}{\mathcal{O}}
\newcommand{\DT}{\operatorname{\mathsf{DT}}}
\newcommand{\DThat}{\operatorname{\widehat{\mathsf{DT}}}}
\newcommand{\HilbBetan}{\operatorname{Hilb}^{\beta ,n}}
\newcommand{\Hilb}{\operatorname{Hilb}}
\newcommand{\Xhat}{\widehat{X}}
\newcommand{\folomo}{{FoLoMo }}

\begin{document}

\begin{abstract}
The topological vertex is a universal series which can be regarded as
an object in combinatorics, representation theory, geometry, or
physics. It encodes the combinatorics of 3D partitions, the action of
vertex operators on Fock space, the Donaldson-Thomas theory of toric
Calabi-Yau threefolds, or the open string partition function of
$\cnums^{3}$.

We prove several identities in which a sum over terms involving the
topological vertex is expressed as a closed formula, often a product
of simple terms, closely related to Fourier expansions of Jacobi
forms. We use purely combinatorial and representation theoretic
methods to prove our formulas, but we discuss applications to the
Donaldson-Thomas invariants of elliptically fibered Calabi-Yau
threefolds at the end of the paper. 
\end{abstract}

\maketitle 

%\markboth{???}  {???}
%\renewcommand{\sectionmark}[1]{}

%\tableofcontents
%\pagebreak

\section{Introduction}\label{sec: intro}

The topological vertex $\Vsf_{\lambda \mu \nu}=\Vsf_{\lambda \mu
\nu}(p)$ is a universal formal Laurent series in $p$ depending on a
triple of partitions $(\lambda, \mu, \nu )$. It can be considered as
an object in combinatorics, representation theory, geometry, or
physics. In combinatorics, $\Vsf_{\lambda \mu \nu}$ is the generating
function for the number of 3D partitions with asymptotic legs of type
$(\lambda, \mu, \nu )$ (see Definition~\ref{defn: box counting
vertex}). In representation theory, $\Vsf_{\lambda \mu \nu}$ is given
as the matrix coefficients of a certain vertex operator on Fock space
\cite{Ok-Re-Va}. In geometry, $\Vsf_{\lambda \mu \nu}$ is the basic
building block for computing the Donaldson-Thomas/Gromov-Witten
invariants of toric Calabi-Yau threefolds \cite{MNOP1}. The
topological vertex was first discovered in physics as an open string
partition function in type IIA string theory on $\cnums^{3}$
\cite{AKMV}. An explicit expression for $\Vsf_{\lambda \mu \nu}$ in
terms of Schur functions was given by \cite{Ok-Re-Va} (see
\S~\ref{sec: topo vertex and schur fncs}).

In this paper we prove several ``trace identities'' in which a sum
over certain combinations of the vertex is expressed as a closed
formula, often a product of simple terms. The products are closely
related to the Fourier expansions of Jacobi forms. These identities
are crucial inputs into the computation of the Donaldson-Thomas
partition functions of certain elliptically fibered Calabi-Yau
threefolds in terms of Jacobi forms \cite{Bryan-K3xE,Bryan-Kool,BOPY}.

\section{Definitions and the main result.}\label{sec: defns and
result}

In this section we give the combinatorial definition of the vertex and
we state our main identities.

\begin{defn}\label{defn: 3D partition asympt to (a,b,c)} Let $(\lambda
,\mu ,\nu )$ be a triple of ordinary partitions. A \emph{3D
partition $\pi $ asymptotic to $(\lambda ,\mu ,\nu )$} is a subset
\[
\pi \subset \left(\znums _{\geq 0} \right)^{3}
\]
satisfying
\begin{enumerate}
\item if any of $(i+1,j,k)$, $(i,j+1,k)$, and $(i,j,k+1)$ is in $\pi
$, then $(i,j,k)$ is also in $\pi $, and
\item
\begin{enumerate}
\item $(j,k)\in \lambda $ if and only if $(i,j,k)\in \pi $ for all $i\gg 0$,
\item $(k,i)\in \mu  $ if and only if $(i,j,k)\in \pi $ for all $j\gg 0$,
\item $(i,j)\in \nu  $ if and only if $(i,j,k)\in \pi $ for all $k\gg 0$.
\end{enumerate}
\end{enumerate}
where we regard ordinary partitions as finite subsets of $\left(\znums
_{\geq 0} \right)^{2}$ via their diagram.
\end{defn}

Intuitively, $\pi $ is a pile of boxes in the positive octant of
3-space.  Condition (1) means that the boxes are stacked stably with
gravity pulling them in the $(-1,-1,-1)$ direction; condition (2)
means that the pile of boxes is infinite along the coordinate axes
with cross-sections asymptotically given by $\lambda $, $\mu $, and
$\nu $.

The subset $\{(i,j,k ): (j,k)\in \lambda \}\subset \pi $ will be
called the \emph{leg} of $\pi $ in the $i$ direction, and the legs in
the $j$ and $k$ directions are defined analogously. Let
\begin{equation*}
\xi _{\pi } (i,j,k) = 1 - \# \text{ of legs of $\pi $ containing }
(i,j,k) .
\end{equation*}

We define the renormalized volume of $\pi $ by
\[
|\pi | = \sum _{(i,j,k)\in \pi } \xi _{\pi } (i,j,k).
\]
Note that $|\pi |$ can be negative.
\begin{defn}\label{defn: box counting vertex}
The topological vertex $\Vsf_{\lambda \mu \nu }$ is defined to be
\[
\Vsf _{\lambda \mu \nu }= \sum _{\pi } p^{|\pi |}
\]
where the sum is taken over all 3D partitions $\pi $ asymptotic to
$(\lambda ,\mu ,\nu )$. We regard $\Vsf _{\lambda \mu \nu }$ as a
formal Laurent series in $p$. Note that $\Vsf _{\lambda \mu \nu }$ is
clearly cyclically symmetric in the indices, and reflection about the
$i=j$ plane yields
\[
\Vsf _{\lambda \mu \nu } = \Vsf _{\mu '\lambda '\nu '}
\]
where $'$ denotes conjugate partition:
\[
\lambda ' = \{(i,j): (j,i)\in \lambda  \}.
\]

\end{defn}
This definition of topological vertex differs from the vertex $C
(\lambda ,\mu ,\nu )$ of the physics literature by a normalization
factor (and we use the variable $p$ instead of $q$). Our $\Vsf
_{\lambda \mu \nu }$ is equal to $P (\lambda ,\mu ,\nu )$ defined by
Okounkov, Reshetikhin, and Vafa \cite[eqn~3.16]{Ok-Re-Va}. They derive
an explicit formula for $\Vsf _{\lambda \mu \nu }=P (\lambda, \mu, \nu
)$ in terms of Schur functions \cite[eqns~3.20 and 3.21]{Ok-Re-Va}.

The \emph{rows} or \emph{parts} of $\lambda $
are the integers $\lambda _j = \min \{i \;|\;(i,j) \not \in \lambda
\}$, for $j \geq 0$. We use the notation
\[
|\lambda | = \sum_{j} \lambda_{j},\quad \| \lambda \| ^{2} =\sum_{j}\lambda_{j}^{2}.
\]
Let $\bx$ denote the partition with a single part of size 1.

We also use the notation
\[
M(p,q) = \prod_{m=1}^{\infty} (1-p^{m}q)^{-m}
\]
and the shorthand $M(p)=M(p,1)$.  Here $M$ stands for MacMahon, who
proved~\cite{macmahon} that 
\[
\Vsf_{\emptyset \emptyset \emptyset} =M(p).
\]

We can now state our main result.

\begin{theorem}\label{thm: main formulas}
The following identities hold as formal power series in $q$ whose
coefficients are formal Laurent series in $p$:
\begin{align}
\sum_{\lambda} q^{|\lambda |}& p^{\| \lambda' \| ^{2}} \Vsf_{\lambda'
\lambda \emptyset }=  M(p) \prod_{d=1}^{\infty} (1-q^{d})^{-1}M(p,q^{d})\label{eqn 2}\\
\quad\nonumber \\
\sum_{\lambda} q^{|\lambda |}&\frac{\Vsf_{\lambda
\bx\emptyset}}{\Vsf_{\lambda \emptyset \emptyset}} =
(1-p)^{-1}\prod_{d=1}^{\infty} \frac{(1-q^{d})}{(1-pq^{d})(1-p^{-1}q^{d})}\label{eqn 3}\\
\quad\nonumber \\
\sum_{\lambda} q^{|\lambda |}& p \, \frac{\Vsf_{\bx \bx 
\lambda}}{\Vsf_{\emptyset \emptyset \lambda}} =
\prod_{m=1}^{\infty}(1-q^{m})^{-1}\cdot \left\{
1+\frac{p}{(1-p)^{2}}+\sum_{d=1}^{\infty}\sum_{k|d}k(p^{k}+p^{-k})q^{d}\right\}\label{eqn 4}\\
\quad&\quad \nonumber \\
\sum_{\lambda} q^{|\lambda |}& p^{\| \lambda \| ^{2}}
\Vsf_{\lambda \lambda' \emptyset} \,\,\frac{\Vsf_{\lambda \bx \emptyset}}{\Vsf_{\lambda \emptyset
\emptyset}} =
(1-p)^{-1}M(p)\prod_{d=1}^{\infty}
\frac{M(p,q^{d})}{(1-pq^{d})(1-p^{-1}q^{d})}. \label{eqn 5}
\end{align}
The sums in the left hand sides of the above formulas run over all
partitions. 
\end{theorem}

%We call these formulas ``trace formulas'' since the left hand side can
%be expressed as the traces of certain operators on Fock space. 
%This
%will be made explicit in section~\ref{sec: vertex ops and the pf of
%eqn 5}.

% We note that formula~\eqref{eqn 1} is elementary and well known.
We prove Formula~\eqref{eqn 2} in section~\ref{sec: topo vertex and
schur fncs} using the orthogonality properties of skew Schur
functions. Formulas~\eqref{eqn 3} and \eqref{eqn 4} are proved in
section~\ref{sec: applications of the Bloch-Okounkov thm} using a
theorem of Bloch-Okounkov \cite{Bloch-Okounkov}. The most difficult
identity to prove is equation~\eqref{eqn 5} which we do in
section~\ref{sec: vertex ops and the pf of eqn 5}. There we prove that
the left hand side of equation~\eqref{eqn 5} is given as the trace of
a certain product of operators on Fock space (hence the term ``trace identities'' in the title). To compute the trace, we
use a trick which involves an ``infinite number'' of permutations of
the operators.

\section{The topological vertex and Schur functions}\label{sec: topo
vertex and schur fncs}

Okounkov-Reshetikhin-Vafa derived a formula for the topological vertex
in terms of skew Schur functions. Translating their formulas
\cite[3.20\& 3.21]{Ok-Re-Va} into our notation, we get:
\begin{equation}\label{eqn: ORV formula for vertex}
\Vsf_{\lambda \mu \nu}(p) = M(p) p^{-\half (\| \lambda \| ^{2}+\| \mu'
\| ^{2}+\| \nu \| ^{2})} s_{\nu '}(\ptotheminusrho ) \sum_{\eta} s_{\lambda
'/\eta}(p^{-\nu -\rho})s_{\mu /\eta}(p^{-\nu '-\rho} ).
\end{equation}

Here, $s_{\alpha /\beta}(x_{1},x_{2},\dots )$ is the skew Schur
function (see for example \cite[\S~5]{MacDonald}) and 
\[
\rho =\left(-\half ,-\frac{3}{2},-\frac{5}{2},\dots  \right)
\]
so that $p^{-\nu -\rho}$ is notation for the variable list
\[
p^{-\nu -\rho} = \left(p^{-\nu_{1} +\half },p^{-\nu_{2} +\frac{3}{2}},\dots   \right).
\]

\begin{remark}\label{rem: vertex/M(p) is a rational function in p}
Note that equation~\eqref{eqn: ORV formula for vertex} implies that
$\Vsf_{\lambda \mu \nu}(p)/M(p)$ is a \emph{rational} function in
$p$.\footnote{By \cite[I.3ex2]{MacDonald}, $s_{\lambda}(p^{\rho})$ is
a rational function in $p^{\half}$. Then since the variable list $p^{-\nu
-\rho}$ differs from the variable list $p^{-\rho}$ in a finite number
of spots, we conclude that $s_{\lambda}(p^{-\nu -\rho})$ is also a
rational function. It follows that $s_{\lambda /\mu}(p^{-\nu -\rho})$
is also a rational function.} If we divide both sides of
equations~\eqref{eqn 2} and \eqref{eqn 5} by $M(p)$, we can regard the
main theorem as identities of power series in $q$ whose coefficients
are rational functions in $p$.
\end{remark}

We prove equation~\eqref{eqn 2} as follows. 
Using equation~\eqref{eqn: ORV formula for vertex} we see 
\[
\Vsf_{\lambda '\lambda \emptyset} = M(p)p^{-\| \lambda' \| ^{2}}
\sum_{\eta} s_{\lambda /\eta}(\ptotheminusrho )^{2}
\]
and so (using orthogonality of skew Schur functions \cite[28(a) pg
94]{MacDonald} in the second line below) we see
\begin{align*}
\sum_{\lambda} q^{|\lambda |} p^{\| \lambda' \| ^{2}} \Vsf_{\lambda
'\lambda \emptyset} &= M(p)\sum_{\lambda ,\eta} q^{|\lambda |} (s_{\lambda /\eta}(p^{\half },p^{\frac{3}{2}},\dots ))^{2}\\
&=M(p)\prod_{d=1}^{\infty} \left((1-q^{d})^{-1}\prod_{i,j=1}^{\infty}(1-q^{d}p^{i-\half +j-\half })^{-1} \right)\\
&= M(p) \prod_{d=1}^{\infty} (1-q^{d})^{-1}\prod_{m=1}^{\infty}(1-q^{d}p^{m})^{-m}\\
&=M(p)\prod_{d=1}^{\infty} (1-q^{d})^{-1}M(p,q^{d}).
\end{align*}

We also use equation~\eqref{eqn: ORV formula for vertex} to derive the
following key formulas:
\begin{lemma}\label{lem: eqns for Vlambdaboxempty/Vlambdaemptyempty
and Vlambdaboxbox/Vlambdaemptyempty} 
The following hold:
\begin{align*}
p^{\half }\,\,  \frac{\Vsf_{\lambda \bx \emptyset}}{\Vsf_{\lambda
\emptyset \emptyset}} &= \sum_{i=1}^{\infty} p^{-\lambda_{i}+i-\half }\\
p\,\,  \frac{\Vsf_{\lambda \bx \bx}}{\Vsf_{\lambda
\emptyset \emptyset}} &= 1-\left(\sum_{i=1}^{\infty}
p^{-\lambda_{i}+i-\half } \right) \left(\sum_{j=1}^{\infty} p^{\lambda_{j}-j+\half } \right).
\end{align*}
The terms in the parentheses on the right hand side of the second
equation should be regarded as Laurent expansions of rational
functions in $p^{\half}$ and multiplication is done as rational
functions (see Remark~\ref{rem: vertex/M(p) is a rational function in
p}). 
\end{lemma}
\proof Applying equation~\eqref{eqn: ORV formula for vertex} to
$\Vsf_{\lambda \bx \emptyset}/\Vsf_{\lambda \emptyset
\emptyset}=\Vsf_{\bx \emptyset \lambda}/\Vsf_{\emptyset \emptyset
\lambda}$ we see that
\[
p^{\half} \frac{\Vsf_{\lambda \bx \emptyset}}{V_{\lambda \emptyset
\emptyset  }} = s_{\bx}(p^{-\lambda -\rho})=s_{\bx}
(p^{-\lambda_{1}+\half },p^{-\lambda_{2}+\frac{3}{2}},\dots ) =
\sum_{i=1}^{\infty} p^{-\lambda_{i}+i-\half }.
\]
Similarly, 
\begin{align*}
p\frac{\Vsf_{\lambda \bx \bx}}{\Vsf_{\lambda \emptyset \emptyset}} =
p\frac{\Vsf_{\bx \bx \lambda}}{\Vsf_{\emptyset \emptyset \lambda }} &=
\sum_{\eta} s_{\bx /\eta}(p^{-\lambda -\rho} )  s_{\bx /\eta}(p^{-\lambda' -\rho} )  \\
&= 1 + s_{\bx}(p^{-\lambda -\rho })s_{\bx}(p^{-\lambda' -\rho }).
\end{align*}

In general we have the following relation (see
\cite[Eqn~(3.10)]{Ok-Re-Va})\footnote{There is a typo in equation~3.10 in
\cite{Ok-Re-Va} --- the exponent on the right hand side should be $-\nu'-\rho$. }
\[
s_{\lambda /\mu}(p^{\nu +\rho}) = (-1)^{|\lambda |-|\mu |}s_{\lambda
'/\mu '}(p^{-\nu '-\rho}).
\]
This equality makes sense as an equality of rational functions in
$p^{\half}$ (see Remark~\ref{rem: vertex/M(p) is a rational function in p}).

In particular 
\begin{equation}\label{eqn: s(p^{nu+rho})=-s(p^{-nu'-rho})}
s_{\bx}(p^{\nu +\rho})=-s_{\bx}(p^{-\nu '-\rho})
\end{equation}
and thus
\begin{align*}
p\frac{\Vsf_{\lambda \bx \bx}}{\Vsf_{\lambda \emptyset \emptyset}} &=
1 - s_{\bx}(p^{-\lambda -\rho })s_{\bx}(p^{\lambda +\rho })\\
&= 1-\left(\sum_{i=1}^{\infty} p^{-\lambda_{i}+i-\half }
\right)\left(\sum_{j=1}^{\infty} p^{\lambda_{j}-j+\half } \right)
\end{align*}
which proves the lemma.

\section{Applications of a theorem of Bloch-Okounkov}
\label{sec: applications of the Bloch-Okounkov thm}

We summarize a result of Bloch-Okounkov \cite{Bloch-Okounkov} and use
it to prove equations~\eqref{eqn 3} and \eqref{eqn 4}.

We define the following theta function
\[
\Theta (p,q) = \eta (q)^{-3}\sum_{n\in \znums}
(-1)^{n}q^{\half (n+\half )^{2}} p^{n+\half } 
\]
where
\[
\eta (q) = q^{\frac{1}{24}}\prod_{m=1}^{\infty}(1-q^{m}).
\]

By the Jacobi triple product formula, $\Theta $ is given by
\[
\Theta (p,q) = (p^{\half} -p^{-\half})\prod_{m=1}^{\infty}
\frac{(1-pq^{m})(1-p^{-1}q^{m})}{(1-q^{m})^{2}}. 
\]

We suppress the $q$ from the notation: $\Theta (p) = \Theta (p,q)$, and
we note that
\[
\Theta (p) = -\Theta (p^{-1}).
\]

\begin{theorem}[Bloch-Okounkov \cite{Bloch-Okounkov}]\label{thm: Bloch-Okounkov thm}
Define the $n$ point correlation function by the formula
\[
F(p_{1},\dots ,p_{n}) = \prod_{m=1}^{\infty}(1-q^{m})\, \sum_{\lambda}
q^{|\lambda|} \prod_{k=1}^{n} \left(\sum_{i=1}^{\infty} p_{k}^{\lambda_{i}-i+\half } \right).
\]
Then
\[
F(p) = \frac{1}{\Theta (p)}
\]
and
\[
F(p_{1},p_{2}) = \frac{1}{\Theta
(p_{1}p_{2})}\left(p_{1}\frac{d}{dp_{1}}\log(\Theta (p_{1}))+ p_{2}\frac{d}{dp_{2}}\log(\Theta (p_{2})) \right).
\]
\end{theorem}

In \cite{Bloch-Okounkov}, formulas for the general $n$ variable
function are given, but we will only need the cases of $n=1$ and
$n=2$.

Using this theorem, we will prove equations~\eqref{eqn 3} and
\eqref{eqn 4} of the main theorem.

\subsection{Proofs of equations~\eqref{eqn 3} and \eqref{eqn
4}}\label{subsec: pfs of eqn 3 and 4} $\, $

We apply Lemma~\ref{lem: eqns for Vlambdaboxempty/Vlambdaemptyempty
and Vlambdaboxbox/Vlambdaemptyempty} and Theorem~\ref{thm: Bloch-Okounkov thm}:
\begin{align*}
\sum_{\lambda}(1-p) q^{|\lambda |} \frac{\Vsf_{\lambda \bx
\emptyset}}{\Vsf_{\lambda \emptyset \emptyset}} &=
(p^{-\half}-p^{\half}) \sum_{\lambda} q^{|\lambda |}\sum_{i=1}^{\infty}p^{-\lambda_{i}+i-\half}\\
&=(p^{-\half}-p^{\half}) \prod_{m=1}^{\infty}(1-q^{m})^{-1} F(p^{-1})\\
&=(p^{-\half}-p^{\half}) \prod_{m=1}^{\infty}(1-q^{m})^{-1} \frac{1}{-\Theta (p)}\\
&=\prod_{m=1}^{\infty} \frac{(1-q^{m})}{(1-pq^{m})(1-p^{-1}q^{m})}
\end{align*}
which proves equation~\eqref{eqn 3}.

Again we apply Lemma~\ref{lem: eqns for
Vlambdaboxempty/Vlambdaemptyempty and Vlambdaboxbox/Vlambdaemptyempty}
and Theorem~\ref{thm: Bloch-Okounkov thm}:
\begin{align*}
\sum_{\lambda} q^{|\lambda |}p \frac{\Vsf_{\lambda \bx
\bx}}{\Vsf_{\lambda \emptyset \emptyset}} &= \sum_{\lambda}
q^{|\lambda
|}\left\{1-\left(\sum_{i=1}^{\infty}p^{-\lambda_{i}+i-\half}
\right)\left(\sum_{j=1}^{\infty}p^{\lambda_{j}-j+\half} \right)  \right\}\\
&= \prod_{m=1}^{\infty}(1-q^{m})^{-1} \left(1-F(p,p^{-1}) \right).
\end{align*}

From Theorem~\ref{thm: Bloch-Okounkov thm}, we see that 
\[
F(p,p^{-1}) = \lim_{ (p_{1},p_{2})\to (p,p^{-1}) } \frac{1}{\Theta
(p_{1}p_{2})} \left(p_{1}\frac{d}{dp_{1}}\log(\Theta (p_{1}))+ p_{2}\frac{d}{dp_{2}}\log(\Theta (p_{2})) \right).
\]
To evaluate this limit, we simplify the above expression. A short
computation shows that 
\[
p\frac{d}{dp}\log(\Theta (p)) = \half \frac{p+1}{p-1}
+\sum_{m=1}^{\infty}\sum_{k=1}^{\infty} \left(-p^{k}+p^{-k} \right)
q^{mk}.
\]
Thus
\begin{align*}
F(p,p^{-1}) =& \lim_{\begin{smallmatrix} (p_{1},p_{2})\to \\
(p,p^{-1}) \end{smallmatrix}}
\left((p_{1}p_{2})^{\half} - (p_{1}p_{2})^{-\half} \right)^{-1}
\prod_{m=1}^{\infty}
\frac{(1-q^{m})^{2}}{(1-(p_{1}p_{2})q^{m})(1-(p_{1}p_{2})^{-1}q^{m})} \\
&\quad \quad \quad \quad \cdot \left\{\half \cdot \frac{p_{1}+1}{p_{1}-1}+\half \cdot \frac{p_{2}+1}{p_{2}-1}
+\sum_{m=1}^{\infty}\sum_{k=1}^{\infty}
\left(-p_{1}^{k}-p_{2}^{k}+p_{1}^{-k}+p_{2}^{-k} \right)q^{mk}
\right\} \\
=&  \lim_{\begin{smallmatrix} (p_{1},p_{2})\to \\
(p,p^{-1}) \end{smallmatrix}}
\frac{-(p_{1}p_{2})^{\half}}{1-p_{1}p_{2}} \cdot
\left\{\frac{p_{1}p_{2}-1}{(p_{1}-1)(p_{2}-1)} +
\sum_{m=1}^{\infty}\sum_{k=1}^{\infty} (1-p_{1}^{k}p_{2}^{k})(p_{1}^{-k}+p_{2}^{-k})q^{mk} \right\}\\
= &  \lim_{\begin{smallmatrix} (p_{1},p_{2})\to \\
(p,p^{-1}) \end{smallmatrix}}
\,(p_{1}p_{2})^{\half} \left\{\frac{1}{(1-p_{1})(1-p_{2})} -
\sum_{m=1}^{\infty }\sum_{k=1}^{\infty } \frac{1-(p_{1}p_{2})^{k}}{1-p_{1}p_{2}}(p_{1}^{-k}+p_{2}^{-k})q^{mk} \right\}\\
=& \frac{1}{(1-p)(1-p^{-1})} - \sum_{m=1}^{\infty}\sum_{k=1}^{\infty}
k(p^{k}+p^{-k})q^{mk}. 
\end{align*}
Therefore 
\[
1-F(p,p^{-1}) = 1 +\frac{p}{(1-p)^{2}} + \sum_{d=1}^{\infty}\sum_{k|d}k(p^{k}+p^{-k})q^{d}
\]
which finishes the proof of equation~\eqref{eqn 4}.

\section{Vertex operators and the proof of equation~\eqref{eqn
5}}\label{sec: vertex ops and the pf of eqn 5}

There are several sources for vertex operators and the infinite wedge
formalism. For consistency, we will follow the notation and conventions
of \cite[Appendix~A]{Okounkov-InfWedge}.

Let $V$ be the vector space with basis $\left\{\underline{k}
\right\},$ $k\in \ZplusHalf$. We define \emph{Fock space} $\FockSpace
$ to be the vector space spanned by vectors
\[
v_{S} = \underline{s_{1}}\wedge \underline{s_{2}}\wedge \dots 
\]
where $S=\left\{s_{1}>s_{2}>\dots \right\}\subset \ZplusHalf$ is
any subset such that the sets
\[
S_{+} = S\cap \left(\ZplusHalf \right)_{>0}\quad \text{and} \quad S_{-} =
S^{c}\cap \left(\ZplusHalf \right)_{<0}
\]
are both finite. Let $(\cdot ,\cdot )$ be the inner product on
$\FockSpace$ such that the basis $\left\{v_{S} \right\}$ is
orthonormal. 

For any $k\in \ZplusHalf $ let $\psi_{k}$ be the operator
\[
\psi_{k}(f) = \underline{k}\wedge f
\]
and let $\psi^{*}_{k}$ be its adjoint.

For any partition $\lambda =\{\lambda_{1}\geq \lambda_{2}\geq \dots
\}$, we define the vector
\[
v_{\lambda}  = \underline{(\lambda_{1}-\half )}\wedge
\underline{(\lambda_{2}-\frac{3}{2})}\wedge \dots  
\]

Let $\FockSpaceZero \subset \FockSpace$ be the subspace spanned by the
vectors $\{v_{\lambda} \}$ where $\lambda$ runs over all
partitions. We call this \emph{charge zero Fock space.}

The \emph{energy operator}
\[
H=\sum_{k>0} k\left(\psi_{k}\psi_{k}^{*}+\psi_{-k}^{*}\psi_{-k} \right)
\]
acts on the basis $v_{\lambda}$ by
\[
Hv_{\lambda}  = |\lambda |v_{\lambda }
\]
and so the operator $q^{H}$ acts by
\[
q^{H}v_{\lambda} = q^{|\lambda |}v_{\lambda }
\]
where $q$ is a formal parameter.

For $n\in \znums$, $n\neq 0$ define
\[
\alpha_{n}=\sum_{k} \psi_{k-n}\psi^{*}_{k}
\]
and observe that $\alpha^{*}_{n}=\alpha_{-n}$.

Following \cite{Okounkov-InfWedge}, we define the \emph{vertex
operators} $\Gamma_{\pm}(\mathbf{x})$ which are operators on
$\FockSpaceZero$ over the coefficient ring given by symmetric
functions in an infinite set of variables $\mathbf{x} =
(x_{1},x_{2},x_{3},\dots )$. Let $\mathbf{s}=(s_{1},s_{2},\dots )$ with
\[
s_{k} =\frac{1}{k}\sum_{i=1}^{\infty} x_{i}^{k}
\]
be the power sum basis for the ring of symmetric functions and
let
\footnote{In \cite{Okounkov-InfWedge}, the argument of
$\Gamma_{\pm}$ is $\mathbf{s}$, and the dependence on the underlying set of
variables $\mathbf{x}$ is left implicit. We prefer to make
$\mathbf{x}$ the explicit argument.}
\[
\Gamma_{\pm}(\mathbf{x}) =
\exp\left(\sum_{n=1}^{\infty}s_{n}\alpha_{\pm n } \right).
\]
Observe that $\Gamma^{*}_{\pm}=\Gamma_{\mp}$.

The matrix coefficients of the vertex operators in the $\{v_{\lambda}
\}$ basis are given by skew Schur functions:
\begin{equation}\label{eqn: matrix coefs of vertex ops are skew schur}
(\Gamma_{-}(\mathbf{x})v_{\mu},v_{\lambda}) =
(v_{\mu},\Gamma_{+}(\mathbf{x})v_{\lambda}) = s_{\lambda
/\mu}(\mathbf{x}).
\end{equation}

Orthogonality of the skew Schur functions then gives rise to the
following commutation equation:
\[
\Gamma_{+}(\mathbf{x})\Gamma_{-}(\mathbf{y}) = \prod_{i,j}
(1-x_{i}y_{j})^{-1} \Gamma_{-}(\mathbf{y})\Gamma_{+}(\mathbf{x}),
\]
in particular
\begin{equation}\label{eqn: Gamma+Gamma- commutation relation}
\Gamma_{+}(u\ptotheminusrho )\Gamma_{-}(v\ptotheminusrho ) =M(p,uv)
\Gamma_{-}(v\ptotheminusrho )\Gamma_{+}(u\ptotheminusrho )
\end{equation}
where recall that $u\ptotheminusrho  = (up^{\frac{1}{2}},up^{\frac{3}{2}},up^{\frac{5}{2}},\dots )$.

We let
\[
\psi (z) = \sum_{i} z^{i}\psi_{i}\quad \text{and}\quad \psi^{*}(w) =
\sum_{j} w^{-j} \psi_{j}^{*}. 
\]

The commutation relations of these operators with the vertex operators
are given by
\begin{align}\label{eqn: commutation of Gamma(x) with psi(z) and psi*(w)}
\Gamma_{\pm}(\mathbf{x})\psi (z) &= \prod_{i=1}^{\infty}
(1-x_{i}z^{\pm 1})^{-1} \,\, \psi (z)\Gamma_{\pm}(\mathbf{x})\\
\Gamma_{\pm}(\mathbf{x})\psi^{*} (w) &= \prod_{i=1}^{\infty}
(1-x_{i}w^{\pm 1}) \,\,\psi^{*}(w)\Gamma_{\pm}(\mathbf{x}).\nonumber
\end{align}

We use operators $\E_{r}$ introduced by Okounkov-Pandharipande in
\cite[\S~2.2.4]{Okounkov-Pandharipande-completed-cycles}.  For $r\in
\znums$, let\footnote{We avoid the use of normal ordering by allowing
coefficients which are Laurent series in $p^{\half}$.}
\[
\E_{r}(p) = \sum_{k\in \ZplusHalf} p^{k-\frac{r}{2}}\,\, \psi_{k-r}\psi^{*}_{k}.
\]
Our variable $p$ is related to the variable $z$ in
\cite{Okounkov-Pandharipande-completed-cycles} by $p=e^{z}.$

From \cite[Eqns~2.9 and 0.18]{Okounkov-Pandharipande-completed-cycles}, we
see that $\E_{0}$ is the diagonal operator given by
\[
\E _0 (p) v_{\lambda} = \left(\sum_{i=1}^{\infty}
p^{\lambda_{i}-i+\half} \right)v_{\lambda}. 
\]
By equation~\eqref{eqn: s(p^{nu+rho})=-s(p^{-nu'-rho})} we have
\begin{equation}\label{eqn: formula for the operator E0}
\E _0 (p) v_{\lambda} = -\left(\sum_{i=1}^{\infty}
p^{-\lambda'_{i}+i-\half} \right)v_{\lambda}. 
\end{equation}

We define
\[
\E (a,p) = \sum_{r\in \znums} a^{r} \E_{r}(p)
\]
where $a$ is a formal parameter. A short computation shows that 
\[
\E (a,p) = \psi (a^{-1}p^{\half})\psi^{*} (a^{-1}p^{-\half}).
\]
From equation~\eqref{eqn: commutation of Gamma(x) with psi(z) and
psi*(w)} we get
\[
\Gamma_{\pm}(\mathbf{x})\E (a,p) =\prod_{i=1}^{\infty}\frac{(1-a^{\mp
1}p^{\mp \half}x_{i})}{(1-a^{\mp
1}p^{\pm \half}x_{i})}\,\, \E (a,p)\Gamma_{\pm}(\mathbf{x}). 
\]

For
$\mathbf{x}=u\ptotheminusrho =(up^{\frac{1}{2}},up^{\frac{3}{2}},up^{\frac{5}{2}},\dots )
$ the above simplifies to 
\begin{align}\label{eqn: commutation of Gamma+- with E(a,p)}
\Gamma_{+}(u\ptotheminusrho ) \E (a,p) &= (1-a^{-1}u) \E (a,p) \Gamma_{+}(u\ptotheminusrho )\\
\E (a,p) \Gamma_{-}(u\ptotheminusrho ) &= (1-au) \Gamma_{-}(u\ptotheminusrho )
\E (a,p). \nonumber
\end{align}
Finally, it follows from equation~\eqref{eqn: matrix coefs of vertex
ops are skew schur} that
\begin{equation}\label{eqn: commutation of Gamma and q^{H}}
\Gamma_{\pm}(\mathbf{x})q^{H} = q^{H} \Gamma_{\pm}(q^{\pm 1}\mathbf{x}). 
\end{equation}

We now write the left hand side of equation~\eqref{eqn 5} in the main
theorem as a trace of operators on charge zero Fock space.
\begin{lemma}\label{lem: eqn 5 written as a trace}
\begin{align*}
\sum_{\lambda} q^{|\lambda |} p^{\| \lambda \| ^{2}} \Vsf_{\lambda
\lambda' \emptyset}\,\, \frac{\Vsf_{\lambda \bx
\emptyset}}{\Vsf_{\lambda \emptyset \emptyset}}& = -p^{-\half} \tr
\left(\E
_0(p)\Gamma_{+}(\ptotheminusrho )\Gamma_{-}(\ptotheminusrho )q^{H} \right)\\
& =- p^{-\half} \operatorname{Coeff}_{a^{0}}\left\{\tr
\left(\E(a,p)\Gamma_{+}(\ptotheminusrho )\Gamma_{-}(\ptotheminusrho
)q^{H} \right) \right\}.
\end{align*}
\end{lemma}
\proof The following computation uses, in order, the commutation
relation for $\Gamma_{+}$ and $\Gamma_{-}$ (equation~\eqref{eqn:
Gamma+Gamma- commutation relation}), the definition of trace, the
formula for $\E_{0} (p)$ (equation~\eqref{eqn: formula for the
operator E0}), Lemma~\ref{lem: eqns for
Vlambdaboxempty/Vlambdaemptyempty and
Vlambdaboxbox/Vlambdaemptyempty}, equation~\eqref{eqn: matrix coefs of
vertex ops are skew schur}, equation~\eqref{eqn: ORV formula for
vertex}, and finally switching $\lambda$ to $\lambda '$ in the sum:
\begin{align*}
\quad &\,\, \quad - p^{-\half} \tr \left(\E_{0}(p)\Gamma_{+}(\ptotheminusrho )\Gamma_{-}(\ptotheminusrho )q^{H} \right)\\
& =- p^{-\half} M(p) \tr \left(\E_{0}(p)\Gamma_{-}(\ptotheminusrho )\Gamma_{+}(\ptotheminusrho )q^{H} \right)\\
&=-p^{-\half}M(p)\sum_{\lambda} \left(v_{\lambda}, \E_{0}(p)\Gamma_{-}(\ptotheminusrho )\Gamma_{+}(\ptotheminusrho )q^{H} v_{\lambda}\right)\\
&=p^{-\half }M(p) \sum_{\lambda} q^{|\lambda |}
\left(\sum_{i=1}^{\infty}
p^{-\lambda'_{i}+i-\half} \right) \left(v_{\lambda}, \Gamma_{-}(\ptotheminusrho )\Gamma_{+}(\ptotheminusrho )v_{\lambda} \right)\\
&=M(p)\sum_{\lambda}q^{|\lambda |} \frac{\Vsf_{\lambda' \bx
\emptyset}}{\Vsf_{\lambda' \emptyset \emptyset}} \left(\Gamma_{+}(\ptotheminusrho )v_{\lambda },\Gamma_{+}(\ptotheminusrho )v_{\lambda } \right)\\
&= \sum_{\lambda} q^{|\lambda |} \frac{\Vsf_{\lambda' \bx
\emptyset}}{\Vsf_{\lambda' \emptyset \emptyset}} M(p)\sum_{\eta} \left(s_{\lambda /\eta}(\ptotheminusrho ) \right)^{2}\\
&= \sum_{\lambda} q^{|\lambda |} \frac{\Vsf_{\lambda' \bx
\emptyset}}{\Vsf_{\lambda' \emptyset \emptyset}}\,  p^{\| \lambda' \| ^{2}}
V_{\lambda '\lambda \emptyset}\\
&= \sum_{\lambda} q^{|\lambda |}  p^{\| \lambda \| ^{2}}
V_{\lambda \lambda' \emptyset}\,  \frac{\Vsf_{\lambda \bx
\emptyset}}{\Vsf_{\lambda \emptyset \emptyset}} \qedhere
\end{align*}

While the operator $\E_{0}(p)$ does not have good commutation
relations with the vertex operators, the operator $\E (a,p)$
does. Hence we first replace $\E_{0}(p)$ with the more general $\E
(a,p)$, compute the trace, and then specialize to the $a^{0}$
coefficient.
\begin{lemma}\label{lem: trace of E Gamma+Gamma-q^{H}}
\begin{multline*}
\tr \left(\E (a,p)\Gamma_{+}(\ptotheminusrho )\Gamma_{-}(\ptotheminusrho )q^{H}
\right) = \\
\frac{M(p)}{p^{\half}-p^{-\half}}\prod_{m=1}^{\infty}
\frac{(1-q^{m}a^{-1})(1-q^{m-1}a)(1-q^{m})
M(p,q^{m})}{(1-pq^{m})(1-p^{-1}q^{m})}.
\end{multline*}
\end{lemma}
\proof Our strategy is the following. We use the cyclic invariance of
trace along with the commutation relations for $\Gamma_{+}$ to move
the operator $\Gamma_{+}$ past the other operators cyclically to the
right until the operators are back to their original positions, but
with new arguments. We perform this operation a countable number of
times, eventually making the $\Gamma_{+}$ operator
disappear\footnote{The third author thanks Guillaume Chapuy and Sylvie
Corteel for teaching him this trick at a conference lunch in
2014. Bouttier, Chapuy, and Corteel used the trick in the paper
\cite{Bouttier-Chapuy-Corteel} in the proof of Theorem 12
therein.}. We then employ the same strategy moving $\Gamma_{-}$
cyclically to the left a countable number of times until it disappears
and we are left with a term which we can evaluate using the
Bloch-Okounkov theorem (Theorem~\ref{thm: Bloch-Okounkov thm}).

We first cyclically commute the operator $\Gamma_{+}$ to the right
using equations~\eqref{eqn: Gamma+Gamma- commutation relation},
\eqref{eqn: commutation of Gamma+- with E(a,p)}, and \eqref{eqn:
commutation of Gamma and q^{H}}:

\begin{align*}
\quad &\quad \quad \quad \,\,\, \tr (\E (a,p)\Gamma_{+}(\ptotheminusrho)\Gamma_{-}(\ptotheminusrho
)q^{H})\\
& = M(p)\tr (\E (a,p)\Gamma_{-}(\ptotheminusrho)\Gamma_{+}(\ptotheminusrho
)q^{H})\\
& = M(p)\tr (\E (a,p)\Gamma_{-}(\ptotheminusrho)q^{H}\Gamma_{+}(q\ptotheminusrho
))\\
& = M(p)\tr (\Gamma_{+}(q\ptotheminusrho
)\E (a,p)\Gamma_{-}(\ptotheminusrho)q^{H})\\
& = M(p)(1-qa^{-1})\tr (\E (a,p)\Gamma_{+}(q\ptotheminusrho) \Gamma_{-}(\ptotheminusrho)q^{H}).
\end{align*}

Cyclically commuting $\Gamma_{+}$ to the right a second time we get:
\begin{multline*}
 \tr (\E (a,p)\Gamma_{+}(\ptotheminusrho)\Gamma_{-}(\ptotheminusrho
)q^{H}) = \\
M(p)(1-qa^{-1})M(p,q)(1-q^{2}a^{-1}) \tr (\E (a,p)\Gamma_{+}(q^{2}\ptotheminusrho)\Gamma_{-}(\ptotheminusrho
)q^{H}).
\end{multline*}

After performing $N$ iterations of this strategy, we arrive at
\begin{multline*}
 \tr (\E (a,p)\Gamma_{+}(\ptotheminusrho)\Gamma_{-}(\ptotheminusrho
)q^{H}) = \\
\prod_{d=1}^{N}M(p,q^{d-1})(1-q^{d}a^{-1}) \tr (\E (a,p)\Gamma_{+}(q^{N}\ptotheminusrho)\Gamma_{-}(\ptotheminusrho
)q^{H}).
\end{multline*}
It follows from equation~\eqref{eqn: matrix coefs of vertex ops are
skew schur} that 
\[
\Gamma_{\pm}(q^{N}\ptotheminusrho )\equiv \operatorname{Id} \mod q^{N}.
\]
So the above two equations imply that the equation
\begin{multline*}
 \tr (\E (a,p)\Gamma_{+}(\ptotheminusrho)\Gamma_{-}(\ptotheminusrho
)q^{H}) = \\
\prod_{d=1}^{\infty }M(p,q^{d-1})(1-q^{d}a^{-1}) \tr (\E (a,p)\Gamma_{-}(\ptotheminusrho
)q^{H})
\end{multline*}
holds to all orders in $q$ and is hence true as a formal power series
in $q$.

We now apply the same strategy commuting $\Gamma_{-}$ to the left:
\begin{align*}
\tr (\E (a,p)\Gamma_{-}(\ptotheminusrho )q^{H}) &=(1-a)\tr (\E (a,p)\Gamma_{-}(q\ptotheminusrho )q^{H}) \\
&=(1-a)(1-aq)\tr (\E (a,p)\Gamma_{-}(q^{2}\ptotheminusrho )q^{H}) \\
&=\dots \\
&=\prod_{d=1}^{\infty}(1-aq^{d-1})\tr (\E (a,p)q^{H})
\end{align*}
and so we have proved
\begin{multline*}
 \tr (\E (a,p)\Gamma_{+}(\ptotheminusrho)\Gamma_{-}(\ptotheminusrho
)q^{H}) = \\
\prod_{d=1}^{\infty }M(p,q^{d-1})(1-q^{d}a^{-1})(1-aq^{d-1})\tr (\E (a,p)q^{H}).
\end{multline*}

From the definition of $\E_{r}(p)$ we see that its matrix entries are
all off-diagonal if $r\neq 0$. Therefore
\begin{align*}
\tr (\E (a,p) q^{H})&=\tr (\E_{0}(p)q^{H})\\
&= \sum_{\lambda} q^{|\lambda |}\sum_{i=1}^{\infty}p^{\lambda_i -i+\half}\\
&=(p^{\half}-p^{-\half})^{-1}\prod_{m=1}^{\infty} \frac{(1-q^{m})}{(1-pq^{m})(1-p^{-1}q^{m})}
\end{align*}
where the last equality follows from the computation in the proof of
equation~\eqref{eqn 3} in \S~\ref{subsec: pfs of eqn 3 and
4}. Combining this with the previous computations finishes the proof
of the lemma. \qed

Combining Lemmas~\ref{lem: eqn 5 written as a trace} and \ref{lem:
trace of E Gamma+Gamma-q^{H}}, we get
\begin{multline*}
\sum_{\lambda} q^{|\lambda |} p^{\| \lambda \|  ^{2}} \Vsf_{\lambda \lambda' \emptyset} \frac{\Vsf_{\lambda
\bx \emptyset}}{\Vsf_{\lambda \emptyset \emptyset}} =\frac{1}{1-p} M(p)\prod_{m=1}^{\infty}
\frac{M(p,q^{m})}{(1-pq^{m})(1-p^{-1}q^{m})}\\
 \cdot \operatorname{Coeff}_{a^{0}}\left\{\prod_{m=1}^{\infty}
(1-q^{m}a^{-1})(1-q^{m-1}a)(1-q^{m}) \right\}.
\end{multline*}

By the Jacobi triple product identity, we have
\[
\prod_{m=1}^{\infty} (1-q^{m}a^{-1})(1-q^{m-1}a)(1-q^{m}) =
\sum_{n=-\infty}^{\infty} q^{\binom{n}{2}} (-a)^{n}
\]
whose $a^{0}$ coefficient is 1. Plugging into the previous equation we
finish the proof of equation~\eqref{eqn 5}.

The proofs of all the formulas in Theorem~\ref{thm: main formulas} are
now complete. 

\section{Geometry and Applications}\label{sec: geometry and applications}

\subsection{Overview.} In this section we briefly outline the
applications of our trace formulas to the Donaldson-Thomas theory of
elliptically fibered Calabi-Yau threefolds
\cite{Bryan-K3xE,Bryan-Kool,BOPY} and the connection with Jacobi
forms. This section is logically independent from the previous
sections and can be safely ignored by readers not interested in the
geometric applications of our theorem.

As we will explain below, each formula in our main theorem gives a
certain contribution to the Donaldson-Thomas partition function of
elliptically fibered threefolds. Namely, we have the following.
\begin{itemize}
\item The contribution of multiples of a nodal elliptic fiber is given
by equation~\eqref{eqn 2}.
\item The contribution of multiples of a smooth elliptic fiber
attached to a smooth section curve is given by equation~\eqref{eqn 3}.
\item The contribution of multiples of a smooth elliptic fiber
attached to a nodal section curve is given by equation~\eqref{eqn 4}.
\item The contribution of multiples of a nodal elliptic fiber
attached to a smooth section curve is given by equation~\eqref{eqn 5}.
\end{itemize}

\subsection{A very short description of Donaldson-Thomas theory.}

Donaldson-Thomas theory is a curve counting theory for Calabi-Yau
threefolds which is conjecturally equivalent to Gromov-Witten theory
\cite{MNOP1}. Let $X$ be a Calabi-Yau threefold, let $\beta \in
H_{2}(X,\znums )$ be a curve class, and let $n\in \znums$.

Let
\[
\HilbBetan (X) = \left\{Z\subset X,\quad [Z]=\beta ,\quad \chi (\O_{Z})=n \right\}
\]
be the Hilbert scheme parameterizing subschemes of class $\beta$ and
whose structure sheaf has holomorphic Euler characteristic $n$. The
Donaldson-Thomas invariant $\DT_{\beta ,n}(X)$ is defined to be the
Behrend function weighted Euler characteristic of the Hilbert scheme:
\[
\DT_{\beta ,n} (X)= \sum_{k\in \znums} k\, e\left( \nu^{-1}(k) \right)
\] 
where $e(-)$ is topological Euler characteristic and 
\[
\nu : \HilbBetan (X)\to \znums 
\]
is Behrend's constructible function \cite{Behrend-micro}. The
unweighted Euler characteristics are also often considered:
\[
\DThat_{\beta ,n} (X) = e\left(\HilbBetan (X) \right).
\]

In this discussion of our applications we will only consider the
unweighted invariants $\DThat_{\beta ,n}(X)$. For further discourse on
including the Behrend function see \cite{Bryan-K3xE,Bryan-Kool,BOPY}.

The invariants are usually assembled into a generating function called
the partition function:
\[
\DThat (X) = \sum_{n,\beta} \DThat_{\beta ,n}(X) q^{\beta} p^{n}.
\]
The above series can be considered as an element in the Novikov
ring of $X$ with coefficients which are formal Laurent series in
$p$.\footnote{With an appropriate choice of a basis of
$H_{2}(X,\znums$), the partition function can be considered as a formal
power series in a set of variables $q_{i}$ whose coefficients are
Laurent series in $p$.}

\subsection{Donaldson-Thomas invariants of a threefold with a torus action}

The topological Euler characteristic of a $\cnums$-scheme $S$ with the
action of a complex torus $T \cong (\cnums^{*})^{k}$ is given by the
Euler characteristic of the fixed locus:
\[
e(S) = e\left(S^{T} \right).
\]
Thus if $X$ is a Calabi-Yau threefold equipped with the action of a
torus $T$, then
\[
\DThat_{\beta ,n}(X) = e\left(\HilbBetan (X)^{T} \right).
\]
In the case where $X$ is toric with the action of a torus $T$, then
$\HilbBetan (X)^{T}$ parameterizes $T$-invariant subschemes and is a
finite set. Thus $\DThat_{\beta ,n}(X)$ simply counts the number of
$T$-invariant subschemes $Z\subset X$ with $[Z]=\beta$ and $\chi
(\O_{Z})=n$.

On a $T$-invariant coordinate chart $\cnums^{3}\subset X$, a
$T$-invariant subscheme is given by a monomial ideal $I\subset \cnums
[x,y,z]$. There is a bijection between 3D partitions and monomial
ideals $\pi \leftrightarrow I$ given by:
\[
(i,j,k)\in \pi \quad \Longleftrightarrow \quad x^{i}y^{j}z^{k} \not \in  I.
\]
If the partition $\pi$ is asymptotic to $(\lambda, \mu, \nu )$ as in
Definition~\ref{defn: 3D partition asympt to (a,b,c)}, then the
subscheme defined by the corresponding monomial ideal $I$ is supported
on the coordinate axes of $\cnums^{3}$ and has nilpotent thickenings
along the axes determined by $(\lambda, \mu, \nu )$.

So we see that in order to compute $\DThat_{\beta ,n}(X)$ for a toric
Calabi-Yau threefold $X$, one needs to count 3D partitions (one for
each coordinate chart) such that the corresponding subschemes agree on
the coordinate overlaps and the curve class and holomorphic Euler
characteristic are given by $\beta $ and $n$ respectively. This leads
to a general formalism developed in \cite{MNOP1} for computing $\DThat
(X)$ in terms of the topological vertex.

\subsection{Using the topological vertex in non-toric geometries}

Even if $X$ has no torus action, under some favorable circumstances,
we can still use the topological vertex technology to compute
Donaldson-Thomas invariants. The idea is to exploit the \emph{motivic}
nature of Euler characteristic. The Euler characteristic is a ring
homomorphism from the Grothendieck group of varieties over $\cnums$ to
the integers. More prosaically, this means that one can compute the
Euler characteristic of a scheme by stratifying it and then summing
the Euler characteristics of each individual stratum.

For example, one can stratify the Hilbert scheme by specifying the support of the
corresponding subschemes. Let $C\subset X$ be a (not necessarily
irreducible) curve, and let
\[
\HilbBetan (X,C) \subset \HilbBetan (X)
\]
be the locus parameterizing subschemes $Z$ such that the reduced
support of $Z$ is contained in $C$. We can then define the
\emph{contribution of $C$ to the Donaldson-Thomas invariants} as
\[
\DThat (X,C) = \sum_{n,\beta} e(\HilbBetan (X,C))p^{n}q^{\beta}.
\]
A key observation is that $\HilbBetan (X,C)$ only depends on
$\Xhat_{C}$, the formal neighborhood of $C$ in $X$. So while $X$
itself may not admit a non-trivial $\cnums^{*}$-action, it is
sometimes the case that $\Xhat_{C}$ does admit an action which then
induces an action on the stratum $\HilbBetan (X,C)$.

This strategy can potentially be iterated: the fixed locus
$\HilbBetan (X,C)^{\cnums^{*}}$ may be further stratified (for example
by the support of the embedded points) and these substrata may admit
further $\cnums^{*}$-actions. The upshot is that in circumstances
where this strategy is successful, it reduces the computation of the
Euler characteristic of the Hilbert scheme to the substratum of the
Hilbert scheme parameterizing subschemes which are \emph{formally
locally monomial ({\folomo } for short)}, i.e. given by monomial
ideals in the formal neighborhoods of each point. Denoting this
substratum by the subscript \folomo, we get
\begin{align*}
\DThat (X,C)&= \sum_{\beta ,n} e\left(\HilbBetan (X,C) \right) p^{n}q^{\beta }\\
&= \sum_{\beta ,n} e\left(\HilbBetan_{\folomo } (X,C) \right) p^{n}q^{\beta }.
\end{align*}
Due to the bijection between monomial ideals and 3D partitions,
\folomo ideals in the formal neighborhood of a point can be counted
with the topological vertex.

\subsection{The case of elliptically fibered Calabi-Yau threefolds.}

Let $X$ be a Calabi-Yau threefold admitting an elliptic fibration $\pi
:X\to S$. We consider the following curves in $X$. Let $F\subset X$ be
some non-singular fiber and let $N\subset X$ be a fiber having a nodal
singularity. Let $B$ be a smooth section curve of genus $g$, that is a
non-singular curve of genus $g$ meeting $F$ and $N$ once each in a
node. The strategy outlined in the previous subsection works
particularly well in this setting where we consider curve classes
which are arbitrary multiples of the fiber class plus possibly a
single multiple of the section class. In these cases, the local
contributions can be expressed in terms of the topological vertex. The
results are the following:

\begin{align*}
\DThat (X,F)&= \sum_{n}\sum_{d=0}^{\infty} e\left(\Hilb_{\folomo}^{d[F],n}(X,F) \right)p^{n}q^{d}\\
&= \sum_{\lambda} q^{|\lambda |}, \\ 
\DThat (X,B+F)&= \sum_{n}\sum_{d=0}^{\infty} e\left(\Hilb_{\folomo}^{[B]+d[F],n}(X,B+F) \right)p^{n}q^{d}\\
&=  \sum_{\lambda} 
q^{|\lambda |} p^{1-g}\, \frac{\Vsf_{\lambda \bx
\emptyset}}{\Vsf_{\lambda \emptyset \emptyset}}\cdot (\Vsf_{\bx \emptyset \emptyset})^{1-2g}, \\
\DThat (X,N)&= \sum_{n}\sum_{d=0}^{\infty} e\left(\Hilb_{\folomo}^{d[N],n}(X,N) \right)p^{n}q^{d}\\
&= \sum_{\lambda} q^{|\lambda |} p^{||\lambda ||^{2}} \Vsf_{\lambda \lambda' \emptyset }, \displaybreak \\
\DThat (X,B+N)&= \sum_{n}\sum_{d=0}^{\infty} e\left(\Hilb_{\folomo}^{[B]+d[N],n}(X,B+N) \right)p^{n}q^{d}\\
&= \sum_{\lambda} q^{|\lambda |} p^{||\lambda ||^{2}+1-g}\, \Vsf_{\lambda
\lambda' \emptyset }\cdot \frac{\Vsf_{\lambda \bx
\emptyset}}{\Vsf_{\lambda \emptyset \emptyset}} \cdot (\Vsf_{\bx \emptyset \emptyset })^{1-2g}.
\end{align*}

There is one additional case to consider.  Let $B'$ be nodal section
curve, that is a curve that has a single nodal singularity at $n\in
B'$ and we assume that $B'$ meets $F$ at the point $n$ so that the
singularity of $B'\cup F$ at $n$ is formally locally that of the
coordinate axes. Let $g'$ be the geometric genus of $B'$. Then
\begin{align*}
\DThat (X,B'+F)&= \sum_{n}\sum_{d=0}^{\infty}
e\left(\Hilb_{\folomo  }^{[B']+d[F],n}(X,B'+F) \right)p^{n}q^{d}\\
&= \sum_{\lambda} q^{|\lambda |}\,p^{1-g'}\, \frac{\Vsf_{\lambda \bx
\bx}}{\Vsf_{\lambda \emptyset \emptyset}}\cdot (\Vsf_{\bx \emptyset \emptyset })^{-2g'}.
\end{align*}

All of the above formulas have a sum over (2D)
partitions. Geometrically, this is because a \folomo subscheme which
is a nilpotent thickening of a curve is determined, away from singular
points and embedded points, by a monomial ideal in the formal
coordinates transverse to the curve, which is in turn determined by a
2D partition. Thus in the above formulas, the term in the sum
corresponding to a partition $\lambda$ counts subschemes which have
thickenings determined by $\lambda$ along the fiber curve ($F$ or
$N$).

Each term in the sum is a combination of topological vertexes which
counts the ways in which embedded points can appear in the
subscheme. For example, in the formula for $\DThat (X,B+N)$ the
$\Vsf_{\lambda \lambda '\emptyset}$ counts the ways of adding embedded
points at the node $n\in N$, the $\Vsf_{\lambda \bx \emptyset}$ counts
ways of adding embedded points at $p=B\cap N$, the $\Vsf_{\lambda
\emptyset \emptyset}$ counts the number of ways of adding embedded
points at some arbitrary point of $N-\{n,p \}$, and the $\Vsf_{\bx
\emptyset \emptyset}$ counts the number of ways of adding embedded
points at some arbitrary point of $B-p$.

The reason that the vertex terms in the above formulas appear with a
(possibly negative) exponent deserves explaining.  In the case of a
toric Calabi-Yau threefold, the number of possible locations for
embedded points is finite --- they only occur at the origins of the
torus invariant coordinate charts. In the case of $\folomo$
subschemes, embedded points can occur at an infinite number of
locations. In the example at hand, the embedded points parameterized
by $\Vsf_{\lambda \emptyset \emptyset}$ can occur at any point of the
curve $N$ except $n$ and $p$ and the embedded points parameterized by
$\Vsf_{\bx \emptyset \emptyset}$ can occur at any point of $B$ except
$p$. The possible locations are parameterized by symmetric products of
$N-\{n,p \}$ and $B-p$ respectively.  Symmetric products are very
amenable to our motivic methods. Indeed, symmetric products endow the
Grothendieck group of varieties with a lambda ring structure which is
compatible with the Euler characteristic homomorphism. Because of the
symmetric product formalism in the Grothendieck group, the
$\Vsf_{\lambda \emptyset \emptyset}$ and the $\Vsf_{\bx \emptyset
\emptyset}$ terms appear with exponents given by the topological Euler
characteristic of $N-\{n,p \}$ and $B-p$ respectively:
\[
(\Vsf_{\lambda \emptyset \emptyset})^{e(N-\{n,p \})}(\Vsf_{\bx \emptyset \emptyset})^{e(B-\{p \})}
=\frac{1}{\Vsf_{\lambda \emptyset \emptyset}}\cdot (\Vsf_{\bx \emptyset \emptyset})^{1-2g} .
\]

\subsection{Applying our trace formulas and the appearance of Jacobi forms.}
By equation~\eqref{eqn: ORV formula for vertex}, we see that 
\[
\Vsf_{\bx \emptyset \emptyset} =M(p) (1-p)^{-1}.
\]

After substituting the above into the equations of the previous
section, we may apply the formulas in our main theorem as well as the
well known formula
\[
\sum_{\lambda} q^{|\lambda |}= \prod_{d=1}^{\infty}(1-q^{d})^{-1}
\]
to the formulas in the previous subsection. We get the following
results:
\begin{align*}
\DThat (X,F)&= \prod_{d=1}^{\infty} (1-q^{d})^{-1}\\
\DThat (X,B+F)&= M(p)^{1-2g} \left(p^{\half}-p^{-\half} \right)^{2g-2} \prod_{d=1}^{\infty}
\frac{(1-q^{d})}{(1-pq^{d})(1-p^{-1}q^{d})}\\
\DThat (X,N) &= M(p)\prod_{d=1}^{\infty} (1-q^{d})^{-1} M(p,q^{d})\\
\DThat (X,B+N) &=  \left(p^{\half}-p^{-\half} \right)^{2g-2} M(p)^{2-2g}
\prod_{d=1}^{\infty} \frac{M(p,q^{d})}{(1-pq^{d})(1-p^{-1}q^{d})}\\
\DThat (X,B'+F) &=  \left(p^{\half}-p^{-\half}
\right)^{2g'}M(p)^{-2g'}\prod_{d=1}^{\infty}(1-q^{d})^{-1}\\
&\quad \cdot 
\left\{1+\frac{p}{(1-p)^{2}}+\sum_{m=1}^{\infty} \sum_{k|m}
k(p^{k}+p^{-k})q^{m} \right\}. 
\end{align*}

The above formulas have a close connection with some well known Jacobi
forms and modular forms. Let $p=\exp\left(2\pi iz \right)$ and
$q=\exp\left(2\pi i\tau \right)$. The Jacobi theta function $\Theta$,
the Dedekind eta function $\eta$, the Weierstrass $\wp$ function, and
the Eisenstein series $G_{2}$ are given as follows:

\begin{align*}
\Theta &= \left(p^{\half}-p^{-\half} \right) \prod_{m=1}^{\infty}
\frac{(1-pq^{m})(1-p^{-1}q^{m})}{(1-q^{m})^{2}}\\
\eta &= q^{\frac{1}{24}}\prod_{m=1}^{\infty}(1-q^{m})\\
\wp &= \frac{1}{12} +\frac{p}{(1-p)^{2}} +\sum_{d=1}^{\infty}
\sum_{k|d} k(p^{k}-2+p^{-k}) q^{d}\\
G_{2} &= -\frac{1}{24} + \sum_{d=1}^{\infty} \sum_{k|d} k q^{d}.
\end{align*}

We may then observe that 
\begin{align*}
\frac{\DThat (X,B+F)}{\DThat (X,F)} &= \left(\frac{M(p)}{p^{\half}-p^{-\half}}
\right)^{1-2g} \frac{1}{\Theta} \\
\frac{\DThat (X,B+N)}{\DThat (X,N)} &= \left(\frac{M(p)}{p^{\half}-p^{-\half}}
\right)^{1-2g} \frac{q^{\frac{1}{24}}}{\Theta \eta } \\
\frac{\DThat (X,B'+F)}{\DThat (X,F)} &= \left(\frac{M(p)}{p^{\half}-p^{-\half}}
\right)^{-2g'} \left\{\wp +2G_{2} +1 \right\}.
\end{align*}

The quotients on the left can be interpreted as the partition function
for \emph{connected} Donaldson-Thomas invariants. There is a general
conjecture due to Huang, Katz, and Klemm \cite{Huang-Katz-Klemm} which
predicts that the Donaldson-Thomas partition function for an
elliptically fibered Calabi-Yau, in classes given by a fixed section
class and multiple fiber classes, is given by a Jacobi form. Several
special cases of this conjecture have been proven using the methods
outlined in this section, namely when $X$ is (1) a local elliptic
surface \cite{Bryan-Kool}, (2) a product of a $K3$ surface and an
elliptic curve \cite{Bryan-K3xE}, and (3) a product of an Abelian
surface and an elliptic curve \cite{BOPY}.

\subsection{Acknowledgments}\label{acknowledgments} We thank Paul
Johnson for showing us (on MathOverFlow) how to use the Bloch-Okounkov
result \cite{Bloch-Okounkov} to prove equation~\eqref{eqn 3} in our
Theorem~\ref{thm: main formulas}. We also thank Guillaume Chapuy and
Sylvie Corteel for showing us their trick of doing an infinite number
of cyclic permutations of vertex operators (see \S~\ref{sec: vertex
ops and the pf of eqn 5}). We thank Georg Oberdieck for providing
feedback on an early draft of this paper.

%MK was supported by a PIMS postdoctoral
%fellowship (CRG Geometry and Physics) while at UBC, and NWO-GQT and
%Marie-Sk{\l}odowska Curie Project 656898 while at Utrecht. JB was
%support by NSERC discovery and accelerator grants. 

\bibliography{/Users/jbryan/jbryan/resources/mainbiblio}
%\bibliography{BKY.bbl}
\bibliographystyle{plain}

\end{document}